\documentclass[11pt,A4paper]{article}
\usepackage{latexsym,amssymb,amsmath,amsfonts,amsthm,mathrsfs}
\usepackage{ifpdf}

\newtheorem{thm}{Theorem}[section]
\def\pf{\noindent{\it Proof.} }
\newcommand{\maj}{\mathrm{maj}}
\newcommand{\Drg}{\mathscr{D}}
\def\qed{\nopagebreak\hfill{\rule{4pt}{7pt}}\medbreak}
\numberwithin{equation}{section}

\begin{document}
\begin{center}
{\Large\bf The Balanced Property of  the $q$-Derangement Numbers and
\\[6pt] the $q$-Catalan Numbers}
\end{center}

\begin{center}
William Y.C. Chen$^1$, David G.L. Wang$^2$, and Larry
X.W. Wang$^3$\\[6pt]
Center for Combinatorics, LPMC-TJKLC\\[6pt]
Nankai University, Tianjin 300071, P.R. China\\[6pt]
$^1$chen@nankai.edu.cn,  $^2$wgl@cfc.nankai.edu.cn,
$^3$wxw@cfc.nankai.edu.cn
\end{center}

\begin{abstract}
Based on B\'ona's condition for the balanced property of the number of of cycles of permutations, we give  a general criterion for the balanced property in terms of the generating function of a
statistic. We show that  the $q$-derangement numbers and the
$q$-Catalan numbers satisfy the
balanced property.
\end{abstract}



%---------------------------------------------------------------------
\section{Introduction}

The notion of the balanced property was
introduced by B\'ona \cite{Bona06}.
Let $\xi$ be a statistic on a set $Q_n$ of combinatorial objects, and let
$\xi(k)$ be the number of objects $\pi$ in $Q_n$ such that
$\xi(\pi)=k$. A statistic $\xi$ over the set $Q_n$ is said to
possess the {\em balanced property} if for any $m\ge2$ and $0\le r\le m-1$,
\begin{equation}\label{bl}
\lim_{n\to\infty}\sum_{k}{\xi(k)\over |Q_n|}={1\over m},
\end{equation}
where the sum ranges over all $k$ congruent to $r$ modulo $m$.
We
assume that $m_n$ is the maximum value of $\xi$.
Then the balanced property (\ref{bl})
can be expressed in terms of a finite sum
\[
\lim_{n\to\infty}
\frac{\xi(r)+\xi(m+r)+\cdots+\xi(pm+r)}{|Q_n|}=\frac{1}{m},
\]
where $p=\lfloor(m_n-r)/m\rfloor$.
In other words, the balanced property means
the asymptotically uniform distribution of a statistic $\xi$ modulo $m$.

B\'ona \cite{Bona06} has shown that the balanced property (\ref{bl})
holds for both the number of
cycles over the set of permutations of $[n]=\{1,2,\ldots,n\}$,
and the number of cycles over the set $\Drg_n$ of derangements of  $[n]$.
Furthermore, he proved that the number of cycles
over $\Drg_n$ with each cycle having
length at least $a$ also satisfies the balanced property.
In a subsequent paper \cite{Bona07}, B\'ona proved the balanced property of
the number of parts over compositions of $n$.

The main objective of this paper is to find
more combinatorial objects that satisfy the balanced property
with respect to certain statistics.
We begin with a general criterion for the balanced property
in terms of the generating functions of the statistics.
This criterion enables us to derive the balanced property of
the major index over derangements of $[n]$ and the major index over Catalan words of length $2n$,
which are counted by the $q$-derangement numbers and the
 $q$-Catalan numbers, respectively.
We also show that the flag major index over  type $B_n$
 derangements satisfies the balanced property.

%---------------------------------------------------------------------
\section{A general criterion}

In this section, we give a formulation of
a general criterion for the balanced property, which can be
used to derive the balanced property of the major index over
derangements of $[n]$ and the major index over Catalan words of length $2n$.
The criterion will be stated in terms of a statistic over
a set $Q_n$.
But it can also be rephrased in terms of the generating
 function of the statistic.
 While our criterion
is based on a general setting, the proof
is essentially the same as the proof of B\'ona for the special
  case concerning the number of cycles of permutations.
  The proof is included for the sake of completeness.

\begin{thm}\label{main}
Let $\xi$ be a statistic on a set $Q_n$ of combinatorial objects.
The balanced property (\ref{bl}) holds if for any $1\le j\le m-1$,
\begin{equation}\label{test}
\lim_{n\to\infty}{f_n\left(\omega^j\right)\over |Q_n|}=0,
\end{equation}
where $\omega=e^{2\pi i/m}$ and
\[
f_n(q)=\sum_{\pi\in Q_n}q^{\xi(\pi)}
\]
is the generating function of the statistic $\xi$.
\end{thm}

\pf Let
\begin{equation}\label{T}
T_n=\sum_{0\le j\le m-1}f_n(\omega^{j})\omega^{-jr}.
\end{equation}
Recall that $\xi(k)$ denotes the number of elements $\pi$ in $Q_n$ such that
$\xi(\pi)=k$. So the generating function $f_n(x)$ can be written as
\begin{equation}\label{eq2}
f_n(q)=\sum_{k}\xi(k)q^k,
\end{equation}
where the sum ranges over all $k$ congruent to $r$ modulo $m$.
Substituting (\ref{eq2}) into (\ref{T}),   we obtain that
\begin{equation}\label{eq3}
T_n=\sum_{k}\xi(k)\sum_{0\le j\le m-1}\omega^{(k-r)j}.
\end{equation}
Since
\begin{equation}\label{omega}
\sum_{0\le j\le m-1}\omega^{(k-r)j}
=\begin{cases}%
m, &\textrm{if }  m|(k-r);\\[5pt]
0, &\textrm{else},
\end{cases}
\end{equation}
the double sum (\ref{eq3}) simplifies to
\[
T_n=m\sum_{k}\xi(k).
\]
Consequently, the balanced property (\ref{bl}) can be recast as
\begin{equation}\label{eq4}
\lim_{n\to\infty}{T_n\over |Q_n|}=1.
\end{equation}
Now, in the expression (\ref{T}) of $T_n$,
the summand $f_n(\omega^{j})\omega^{-jr}$
for $j=0$ is $f_n(1)=|Q_n|$. Thus
\begin{equation}\label{eq6}
T_n=|Q_n|+\sum_{1\le j\le m-1}f_n(\omega^{j})\omega^{-jr}.
\end{equation}
Substituting (\ref{eq6}) into (\ref{eq4}), we conclude that
the the balanced property holds if and only if
\begin{equation}\label{eq5}
\lim_{n\to\infty}\sum_{1\le j\le m-1}\frac{f_n(\omega^{j})}{|Q_n|}\omega^{-jr}=0.
\end{equation}
Evidently, the condition (\ref{test}) implies
(\ref{eq5}) since $\left|\omega^{-jr}\right|=1$. This completes the proof. \qed

%---------------------------------------------------------------------
\section{The $q$-derangement numbers}

In this section, we shall show that the major index
over derangements
of $[n]$ satisfies the balanced property.
Let $S_n$ be the set of permutations of the set $[n]$.
The major index of a permutation $\pi=\pi_1\pi_2\cdots\pi_n\in S_n$
is defined to be
the sum of the indices $i$ such that $\pi_i>\pi_{i+1}$, that is,
\[
\maj(\pi)=\sum_{\pi_i>\pi_{i+1}}i.
\]
Denote by $\maj(k)$ the number of permutations of $[n]$ with major index $k$.

The problem on the balanced property of the major index over permutations has been
considered by
Gordon \cite{Gor63} and Roselle \cite{Ros74}. For any
coprime numbers $k,l\le n$, the number of permutations $\pi$ of $[n]$
such that $\maj(\pi)$ is congruent to $i$ modulo $k$ and $\maj(\pi^{-1})$
is congruent to $j$ modulo $l$ equals $n!/(kl)$, which is independent of $i$ and $j$. To be
more specific,
\begin{equation}\label{GR}
\left|\{\pi\,|\  \maj(\pi)\equiv i\!\!\!\!\pmod k,\
\maj\left(\pi^{-1}\right)\equiv j\!\!\!\!\pmod l\}\right|={n!\over kl}.
\end{equation}
Taking $l=1$, the formula (\ref{GR}) specializes  to the fact that the number of
permutations of $[n]$ with major index $i$ modulo $k$ equals $n!/k$.
In other words,
for any $n\ge m$ and any $0\le i\le m-1$, we have
\begin{equation}\label{eq7}
\sum_{k}{\maj(k)\over n!}={1\over m},
\end{equation}
where the sum ranges over all $k$ congruent $i$ modulo $m$. As noted in
\cite{BMS02}, the relation (\ref{GR}) was
 implicit in Gordon \cite{Gor63} and has been made explicit by Roselle \cite{Ros74}.
 When $l$ divides $n-1$ and $k$ divides $n$,
 a combinatorial proof has been given by Barcelo, Maule and Sundaram \cite{BMS02}.
  A more detailed description of the background on the relation (\ref{GR})  can also
   be found in \cite{BMS02}.
The identity (\ref{eq7}) can be viewed as an exact balanced
property in comparison with the balanced property in the asymptotic sense.
Moreover,
 (\ref{GR})
can be considered as a bivariate version of the exact
balanced property.

Using representations of the symmetric group,
Barcelo and Sundaram \cite[Theorem 2.6]{BS93} have obtained (\ref{eq7}) for the
special case $m=n$. They also gave
a bijective proof in this case.
Recently, Barcelo, Sagan, and Sundaram~\cite{BSS07} gave a combinatorial interpretation
of (\ref{GR}) in the general case by using shuffles of permutations.

We note that our proof of Theorem~(\ref{main})
easily applies to the exact balanced property (\ref{eq7}).
In general, we say a statistic $\xi$
processes the {\em exact balanced property} if for any $0\le r\le m-1$ and $n\ge N$,
\[
\sum_{k}{\xi(k)\over |Q_n|}={1\over m},
\]
where the sum ranges over all $k$ congruent to $r$ modulo $m$,
and $N$ depends only on $m$.
Inspecting the derivation of the formula (\ref{eq5}), we see that
the exact balanced property holds if and only if
\begin{equation}\label{eq8}
\sum_{1\le j\le m-1}\frac{f_n\left(\omega^{j}\right)}{|Q_n|}\omega^{-jr}=0.
\end{equation}
In the usual notation $[0]_q!=1$,
$[n]_q!=[1]_q[2]_q\cdots [n]_q$,
where $[n]_q=1+q+\cdots + q^{n-1}$ for $n\geq 1$,
the generating function for the major index of
permutations of $[n]$ is known to be
\[
f_n(q)=\sum_{\pi\in S_n}q^{\maj(\pi)}=[n]_q!,
\]
see Andrews \cite{And98}.
Since for $n\ge m$, $f_n(q)$ contains the factor
$[m]_q$,
$f_n\left(\omega^j\right)$ contains the factor $[m]_{\omega^j}$.
It follows that $f_n\left(\omega^j\right)$ equals zero  for any $1\le j\le m-1$.
Therefore the relation (\ref{eq8}) holds,
which implies the exact balanced property (\ref{eq7}).

Moreover, it is not difficult to derive the balanced property
and the exact balanced property of the flag major index
over permutations of type $B_n$, which is introduced by Adin and Roichman~\cite{AR01}.
Denote by $S_n^B$ the set of $B_n$-permutations. Then the generating function for
the flag major index over  $S_n^B$ is given by
\[
f^B_n(q)=\sum_{\pi\in S_n^B}q^{\mathrm{fmaj}(\pi)}=[2]_q[4]_q\cdots[2n]_q,
\]
see Chow~\cite{Chow06}. For $n\ge2m-1$,
$f^B_n(q)$ contains the factor $[2m]_q$.
It follows from  (\ref{omega}) that $f^B_n\left(\omega^j\right)=0$.
This yields (\ref{eq8}), leading to the exact balanced property
of the flag major index over $B_n$-permutations.

A derangement of $[n]$ is a permutation
$\pi_1\pi_2\cdots\pi_n$ of $[n]$ such that $\pi_i\ne i$
for all $1\le i\le n$. The following counting formula of derangements with respect to the major index was given
 by Gessel and published later in \cite{GR93},
\begin{equation}\label{drg}
d_n(q)=\sum_{\pi\in\Drg_n}q^{\maj(\pi)}
=\sum_{0\le k\le n}(-1)^kq^{k\choose2}\prod_{k+1\le t\le n}[t]_q,
\end{equation}
where $[0]_q!=1$. Combinatorial proofs for (\ref{drg}) has been found by Wachs~\cite{Wac89}, and Chen and Xu~\cite{CX08}.

\begin{thm}\label{thm_drg}
The major index over derangements of $[n]$ satisfies the balanced property.
In other words, for any  $0\le r\le m-1$, we have
\[
\lim_{n\to\infty}%
\frac{\maj(r)+\maj(m+r)+\cdots+\maj(pm+r)}{|\Drg_n|}
=\frac{1}{m},
\]
where
\[
p=\left\lfloor\frac{{n\choose2}-r}{m}\right\rfloor.
\]
\end{thm}

\pf  Consider the values of $d_n(q)$ evaluated at $q=\omega^j$, namely,
\begin{equation}\label{drgw}
d_n\left(\omega^j\right)=\sum_{0\le k\le n}(-1)^k\omega^{j{k\choose2}}%
\prod_{k+1\le t\le n}[t]_{\omega^j}.
\end{equation}
Note that among the $m$ consecutive integers $n,n-1,\ldots,n-m+1$, there
exists an integer $a$ which can be divided by $m$. For such a choice of $a$, we have
$[a]_{\omega^j}=0$. So any summand in (\ref{drgw}) containing the factor
$[a]_{\omega^j}$ can be ignored. Consequently, the formula (\ref{drgw})
 reduces to the following form
\begin{equation}\label{dw}
d_n\left(\omega^j\right)%
=\sum_{n-m+2\le k\le n}(-1)^k\omega^{j{k\choose2}}%
\prod_{k+1\le t\le n}[t]_{\omega^j}.
\end{equation}
In order to estimate $[t]_{\omega^j}$, let us compute $[t]_{\omega^j}^2$.
  For any $t\in\{k+1,k+2,\ldots,n\}$, we have
\begin{equation}\label{t}
\left|[t]_{\omega^j}\right|^2%
=\left|\frac{1-\omega^{tj}}{1-\omega^j}\right|^2%
=\frac{1-\cos(2\pi tj/m)}{1-\cos(2\pi j/m)},%
\end{equation}
which is clearly bounded by $2/c$,
where \[ c=\min\{1-\cos(2\pi j/m)\,|\, 1\le j\le m-1\}.\]
It is clear that  $2/c\ge 1$. Hence,
\[
\prod_{k+1\le t\le n}\left|[t]_{\omega^j}\right|^2
\le\left({2\over c}\right)^{n-k}\le\left({2\over c}\right)^{m-2}.
\]
Observe that the above estimate is independent of $j$. It follows from (\ref{dw}) that
\[
\left|d_n\left(\omega^j\right)\right|\le \sum_{n-m+2\le k\le n}\left|(-1)^k\omega^{j{k\choose2}}%
\prod_{k+1\le t\le n}[t]_{\omega^j}\right|
\le (m-1)\left({2\over c}\right)^{(m-2)/2}.
\]
Thus $d_n\left(\omega^j\right)$ is bounded by a constant.
By Theorem~\ref{main}, the major index over
$\Drg_n$ satisfies the balanced property.
This completes the proof. \qed

For the type $B$ case, denote by $\Drg_n^B$ the set of derangements of type $B_n$.
Adin and Roichman \cite{AR01} have shown that the generating function of
the flag major index over $\Drg_n^B$ equals
\[
\sum_{\pi\in\Drg_n^B}q^{\mathrm{fmaj}(\pi)}
=\sum_{0\le k\le n}(-1)^kq^{k(k-1)}[2n]_q[2n-2]_q\cdots[2k+2]_q,
\]
see also Chow \cite{Chow06}, and Adin, Brenti and Roichman \cite{ABR01}.
By an argument very similar to the proof of Theorem \ref{thm_drg},
we can derive the balanced property
of the flag major index over $B_n$-derangements.

%----------------------------------------------------------------------
\section{The $q$-Catalan numbers}

In this section, we shall derive the balanced property of the major index
over Catalan words of length $2n$.
The $q$-Catalan numbers are defined by
\begin{equation}\label{gf_ctl}
C_n(q)=\frac{1}{[n+1]_q}{2n\brack n}_q
=\prod_{2\le t\le n}\frac{[t+n]_q}{[t]_q},
\end{equation}
where
\[
{n\brack k}_q=\frac{[n]_q!}{[k]_q![n-k]_q!},
\]
see, for example, Andrews \cite{And87}, and F\"{u}rlinger and Hofbauer \cite{FH85}.
A combinatorial
interpretation of $C_n(q)$ in term of the major index of the
Catalan words of length $2n$ has been given in \cite{FH85}.
A {\it Catalan word} $w$ of length $2n$ is a sequence consisting of
$n$ $0$'s and $n$ $1$'s such that no prefix contains more
$1$'s than $0$'s. Denote the set of Catalan words of length $2n$ by
$\mathscr{C}_n$. The {\em major index} for a Catalan word $w=w_1w_2\cdots
w_{2n}\in\mathscr{C}_n$ is defined by
\[
\maj(w)=\sum_{w_i>w_{i+1}}i.
\]
Let $\maj(k)$ be the number of Catalan words of
length $2n$ with major index $k$. F\"{u}rlinger and Hofbauer have shown that
\[
C_n(q)=\sum_{w\in \mathscr{C}_n}q^{\maj(w)}.
\]
The number of Catalan words of length $2n$
 is given by the Catalan number
\begin{equation}\label{Ctl}
C_n(1)=\frac{1}{n+1}{2n\choose n}.
\end{equation}

\begin{thm}
The major index over Catalan words of length $2n$ satisfies
the balanced property. In other words,
for any $0\le r\le m-1$, we have
\[
\lim_{n\to\infty}%
\frac{\maj(r)+\maj(m+r)+\cdots+\maj(pm+r)}{C_n(1)}
=\frac{1}{m},
\]
where
\[
p=\left\lfloor\frac{n(2n-1)-r}{m}\right\rfloor.
\]
\end{thm}

\pf By Theorem~\ref{main}, it suffices to show that for any $0\le j\le m-1$,
\[
\lim_{n\to\infty}\frac{C_n\left(\omega^j\right)}{C_n(1)}=0.
\]
Let $1\le j\le m-1$. Suppose that $j/m=u/v$,
where $u$ and $v$ are coprime positive integers with
$2\le v\le m$. Write $n=lv+s$, where $0\le s\le v-1$.
Denote the denominator of (\ref{gf_ctl})
by $D_n(q)$, namely,
\[
D_n(q)=\prod_{2\le t\le n}[t]_q=[n]_q!.
\]

It should not be overlooked that the denominator $D_n(q)$
vanishes when $q$ is set to $\omega^j$.
In fact, since $\omega^{jv}=e^{2\pi ijv/m}=e^{2\pi iu}=1$,
we have $1-\omega^{jt}=0$ for any $v|t$.
More precisely,
$D_n(q)$ contains the factor
\[
F_n(q)=\prod_{1\le k\le l}[kv]_q=[v]_q^l\cdot[l]_{q^v}!,
\]
in which the factor $[v]_q^l$ causes
$D_n(q)$ to vanish when evaluated at $q=\omega^j$.
We proceed to represent $C_n(q)$
as a quotient whose denominator does not vanish for $q=\omega^j$.

If $0\le s\le v-2$,
let \[ A_n(q) =\frac{\prod_{l+1\le k\le 2l}[kv]_q}{F_n(q)}.\]
Canceling the common factor $[v]_q^l$ in the numerator and the denominator,
we can reduce it to the form
\[
A_n(q)={2l\brack l}_{q^v}.
\]
Denote the quotient $C_n(q) / A_n(q)$ by $B_n(q)$. It can be checked that
\begin{equation}
B_n(q) =\left(\prod_{2\le t\le n\atop{t\not\in\{v-s,2v-s,\ldots,lv-s\}}}[t+n]_q\right)
\left(\prod_{2\le t\le n\atop{v\,\nmid\,t}}[t]_q\right)^{-1}.\label{B}
\end{equation}
Clearly, the denominator of $B_n(q)$  does not vanish for $q=\omega^j$.

For the case $s=v-1$, let \[
U_n(q)=\frac{\prod_{l+1\le k\le 2l+1}[kv]_q}{F_n(q)}.
\]
Similarly, $U_n(q)$ can be reduced to the following form
\begin{equation}\label{U}
U_n(q)={2l+1\brack l}_{q^v},
\end{equation}
where the denominator does not vanish for $q=\omega^j$. Moreover,
we see that
 $V_n(q)=C_n(q) / U_n(q)$ has the following representation
\[
V_n(q)=\left(\prod_{2\le t\le
n\atop{t\not\in\{2v-s,\ldots,lv-s, (l+1)v-s\}}}[t+n]_q\right)
\left(\prod_{2\le t\le n\atop{v\,\nmid\,t}}[t]_q\right)^{-1}.
\]
Again,  the denominator of $V_n(q)$ is nonzero when $q=\omega^j$.

We are now ready to give an estimate of $\left|C_n\left(\omega^j\right)\right|$.
First,  consider the case $0\leq s\leq v-2$. Since $\omega^{jv}=1$, we have
\begin{align}
A_n\left(\omega^j\right)&={2l\choose l},\label{A1}\\
B_n\left(\omega^j\right)&=\left(\prod_{2\le t\le n\atop{t\not\in\{v-s,2v-s,\ldots,lv-s\}}}[t+s]_{\omega^j}\right)
\left(\prod_{2\le t\le n\atop{v\,\nmid\,t}}[t]_{\omega^j}\right)^{-1}.\label{B1}
\end{align}
In order to give an estimate for $\left|B_n\left(\omega^j\right)\right|$, it is necessary to reduce the above
expression to a quotient such that the numbers of factors of the numerator and
the denominator are both finite as $n$ tends infinity.
It is easy to verify that $B_n\left(\omega^j\right)=1$ for $s=0$ and $s=1$.
By (\ref{B1}), we have
\[
B_n\left(\omega^j\right)
=\left([v-1]_{\omega^j}!^l\prod_{s+2\le t\le 2s}[t]_{\omega^j}\right)
\left([v-1]_{\omega^j}!^l\prod_{2\le t\le s}[t]_{\omega^j}\right)^{-1}
=\prod_{2\le t\le s}\frac{[t+s]_{\omega^j}}{[t]_{\omega^j}}.
\]
In view of (\ref{t}), we see that
\[
\prod_{2\le t\le s}\frac{\left|[t+s]_{\omega^j}\right|^2}{\left|[t]_{\omega^j}\right|^2}
=\prod_{2\le t\le s}\frac{1-\cos(2\pi (t+s)j/m)}{1-\cos(2\pi tj/m)}
\le\left(\frac{2}{c}\right)^{s-1},
\]
where
\[
c=\min\{1-\cos(2\pi tj/m)\,|\, 1\le j\le m-1,\,2\le t\le m-2\}.
\]
Consequently, $\left|B_n\left(\omega^j\right)\right|$ is bounded by a constant, say,
$c_1$. By (\ref{A1}), for $0\le s\le v-2$, we get
\begin{equation}\label{lim1}
\left|C_n\left(\omega^j\right)\right|\le c_1{2l\choose l}.
\end{equation}

For the case $s=v-1$, substituting $q=\omega^j$ in (\ref{U}), we obtain that
\begin{equation}\label{Uomega}
U_n\left(\omega^j\right)={2l+1\choose l}.
\end{equation}
On the other hand, by an analogous argument to  the case $s\le v-2$,
it can be shown that $\left|V_n\left(\omega^j\right)\right|$ is also bounded by a constant, say $c_2$. It follows from (\ref{Uomega}) that for $s=v-1$,
\begin{equation}\label{lim2}
\left|C_n\left(\omega^j\right)\right|\le c_2{2l+1\choose l}.
\end{equation}

Up to now, we have obtained the  estimates
for the $\left|C_n\left(\omega^j\right)\right|$ in the above two cases, namely,
(\ref{lim1}) for  $0\le s\le v-2$, and (\ref{lim2}) for  $s=v-1$.
By (\ref{Ctl}), we can derive the following general upper bound
\begin{equation}\label{lim}
\frac{\left|C_n\left(\omega^j\right)\right|}{C_n(1)} \le c_3\cdot\frac{(n+1){{2l+1}\choose
l}}{{2n\choose n}},
\end{equation}
where $c_3=\max(c_1,c_2)$.
Based on Stirling's formula,
the central binomial coefficient
can be estimated as follows
\begin{equation}\label{2n}
{2n\choose n}\sim\frac{2^{2n}}{\sqrt{\pi n}},\quad\mbox{ as } n\to\infty.
\end{equation}
Since $l=(n-s)/v$ tends to infinity as $n\to\infty$,
we have
\begin{equation}\label{2l}
{{2l+1}\choose l}\sim2{2l\choose l}\sim \frac{2^{2l+1}}{\sqrt{\pi
l}},\quad\mbox{ as } n\to\infty.
\end{equation}
Combining (\ref{lim}), (\ref{2n}) and (\ref{2l}), we find that
\[
\lim_{n\to\infty}\frac{\left|C_n\left(\omega^j\right)\right|}{C_n(1)}
\le\lim_{n\to\infty}c_3\cdot\sqrt{\frac{n}{l}}\cdot\frac{(n+1)}{2^{2n-2l-1}}=0.
\]
Thus the balanced property follows from Theorem~\ref{main}.
This completes the proof.
\qed

\noindent{\bf Acknowledgments.} This work was supported by the 973
Project, the PCSIRT Project of the Ministry of Education, and the
National Science Foundation of China.

%-------------------------------------------------------------------

\end{document}